\documentclass[leqno,10pt,draft,a4paper]{amsart}
\usepackage{amsmath}
\usepackage{amssymb}
\usepackage{latexsym}

\DeclareMathOperator{\Hom}{Hom}
\DeclareMathOperator{\Ext}{Ext}

\renewcommand{\ge}{\geqslant}
\renewcommand{\le}{\leqslant}

\newcommand{\OO}{\mathcal{O}}
\newcommand{\C}{\mathbb{C}}
\newcommand{\Q}{\mathbb{Q}}
\newcommand{\N}{\mathbb{N}}
\newcommand{\R}{\mathbb{R}}

\newcommand{\Z}{\mathbb{Z}}
\newcommand{\h}{\mathfrak{h}}
\newcommand{\g}{\mathfrak{g}}
\DeclareMathOperator{\glob}{glob}
\newcommand{\frob}{^{\mathrm F}}
\newcommand{\up}{\uparrow}
\newcommand{\rarr}{\rightarrow}

\DeclareMathOperator{\inj}{inj}
\DeclareMathOperator{\proj}{proj}

\DeclareMathOperator{\good}{{\mathcal F}(\nabla)}
\DeclareMathOperator{\doog}{{\mathcal F}(\Delta)}

\DeclareMathOperator{\gfd}{gfd}
\DeclareMathOperator{\wfd}{wfd}
\DeclareMathOperator{\pr}{pr}

\newcommand{\lrimpl}{\Leftrightarrow}

\newcommand{\Mod}{\mathrm{mod}}

\DeclareMathOperator{\Stab}{Stab}
\DeclareMathOperator{\Ind}{Ind}

\DeclareMathOperator{\partn}{\Lambda^+}
\newcommand{\Wedge}{\mbox{$\bigwedge$}}
\newcommand{\les}{long exact sequence }
\newcommand{\ses}{short exact sequence }
\newcommand{\GL}{\mathrm{GL}}
\newcommand{\SL}{\mathrm{SL}}

\newcommand{\half}{\frac{1}{2}}

\newcommand{\bl}{\bar{l}}

\newcommand{\wmu}{w\cdot\mu}
\newcommand{\wslamb}{ws\cdot\lambda}
\newcommand{\wlamb}{w\cdot\lambda}
\newcommand{\vlamb}{v\cdot\lambda}

\newcommand{\alphac}{\alpha\check{\ }\,}

\newcommand{\betac}{\beta\check{\ }\,}

\begin{document}
\theoremstyle{plain}
\newtheorem{thm}{Theorem}[section]
\newtheorem{propn}[thm]{Proposition}
\newtheorem{cor}[thm]{Corollary}
\newtheorem{clm}[thm]{Claim}
\newtheorem{lem}[thm]{Lemma}
\newtheorem{conj}[thm]{Conjecture}
\theoremstyle{definition}
\newtheorem{defn}[thm]{Definition}
\newtheorem{rem}[thm]{Remark}

\title{On the good filtration dimension of Weyl modules 
for a linear algebraic group}
\author{Alison E. Parker}

\address{School of Mathematics and Statistics F07, University of Sydney,
  NSW 2006, Australia}
\email{alisonp@maths.usyd.edu.au}



\thanks{
This research was supported by the Association of Commonwealth
Universities and the British Council.}

\begin{abstract}
Let $G$ be a linear algebraic group over an algebraically closed field
of characteristic $p$ whose corresponding root system is irreducible.
In this paper we calculate the Weyl filtration dimension of the
induced $G$-modules, $\nabla(\lambda)$ and the simple $G$-modules
$L(\lambda)$, for $\lambda$ a regular weight.
We use this to 
calculate some $\Ext$ groups of the form
$\Ext^*\bigl(\nabla(\lambda),\Delta(\mu)\bigr)$,
$\Ext^*\bigl(L(\lambda),L(\mu)\bigr)$,
 and
$\Ext^*\bigl(\nabla(\lambda), \nabla(\mu)\bigr)$,
where $\lambda , \mu$ are regular and $\Delta(\mu)$ is the Weyl
module of highest weight $\mu$.
We then deduce the projective dimensions and injective dimensions for
$L(\lambda)$, $\nabla(\lambda)$ and $\Delta(\lambda)$ for $\lambda$ a
regular weight in associated generalised Schur algebras.
We also deduce the global dimension of the Schur algebras for $\GL_n$,
$S(n,r)$, when $p>n$ and for $S(mp,p)$ with $m$ an integer. 
\end{abstract}
\maketitle

\section*{Introduction}
In this paper we consider the notion of the Weyl filtration dimension
and good filtration dimension of modules for a linear algebraic group.
These concepts were first introduced by Friedlander and
Parshall~\cite{fripar} and may be considered a variation of the notion
of projective dimension and injective dimension respectively.
(The precise definition is given in~\ref{defn:gfd}.)
The Weyl filtration dimension of a module is always at most its
projective dimension. In fact, it is often much less.
In the situation of algebraic groups the Weyl and good filtration
dimensions are always
finite for a finite dimensional module (unlike the projective and
injective dimensions which are usually infinite).
Thus knowing these dimensions give us another tool for calculating
the cohomology of an algebraic group. Indeed we use knowledge of these
dimensions to calculate various $\Ext$ groups for $G$.

We had previously calculated the good filtration dimension of the
irreducible modules for $S(n,r)$, the Schur algebra corresponding to $\GL_n(k)$ 
when $n=2$ and $n=3$ in~\cite{parker1}. We were then able to determine
the global dimension of $S(n,r)$.
The proof in~\cite{parker1} relies heavily on 
the use of filtrations
of the induced modules $\nabla(\lambda)$, $\lambda$ a dominant weight, by
modules of the form $\nabla(\mu)\frob \otimes L(\nu)$.

In this paper we instead use 
the translation functors introduced by Jantzen to
calculate properties of the induced modules 
and the Weyl modules (denoted $\Delta(\lambda)$) for an
algebraic group. 
We first calculate the Weyl filtration dimension (abbreviated wfd), of the induced
modules for regular weights (theorem~\ref{thm:wfdlamb}). We then prove 
$\Ext^i\bigl(\nabla(\lambda),\Delta(\mu)\bigr)\cong k$
when $i=\wfd\bigl(\nabla(\lambda)\bigr)+\wfd\bigl(\nabla(\mu)\bigr)$
and $\lambda, \mu$ regular (theorem~\ref{thm:extiind}).
We can then deduce that
$\Ext^i\bigl(L(\lambda),L(\mu)\bigr)\cong k$
for $i=\wfd\bigl(L(\lambda)\bigr)+\wfd\bigl(L(\mu)\bigr)$
(corollary~\ref{cor:wfdirr}).
These results then enable us to write down the injective and
projective dimensions of $L(\lambda)$, $\nabla(\lambda)$ and
$\Delta(\lambda)$ for $\lambda$ a regular weight in associated
generalised Schur algebras (theorem~\ref{thm:injproj}).

We can deduce the value of the global dimension of
$S(n,r)$ when $p>n$ and $S(p,mp)$ with $m \in \N$
(theorems~\ref{thm:glob1} and \ref{thm:glob2}). This gives us an 
alternative proof for $S(2,r)$
(all $p$) and for $S(3,r)$ with $p\ge 5$.
Some of this work also appears in the author's PhD
thesis~\cite{mythesis}, chapter 6.

In general the global dimension of $S(n,r)$ is still not known. 
Previous values were calculated for $r \le n$
by Totaro~\cite{totaro} (for the
classical case) and Donkin~\cite{donkbk}, section 4.8, (for the quantum case).
The semi-simple Schur algebras (that is the Schur algebras with zero global
dimension) have been determined
in~\cite{dotynak} for the classical case 
and~\cite{erdnak}, theorem (A), for the quantum case.
Conjectured values for the remaining cases are presented in 
 \cite{mythesis}, section 6.5.

We conclude by showing that analogous results for 
the Dipper--Donkin quantum group hold and hence for the $q$-Schur
algebra.
The extent to which similar methods
may be applied to category $\OO$ is also discussed.

The author thanks her PhD supervisor, Stephen Donkin for his great help and
encouragement as well as Anton Cox and Karin Erdmann for various
comments on preliminary versions of this paper.

\section{Preliminaries}\label{sect:prelim}
We first review the basic concepts and most of the notation that we will be using.
The reader is referred to
\cite{humph} and \cite{springer} for further information. 
This material is also in~\cite{jantz} where is it presented in the form of group
schemes.

Throughout this paper $k$ will be an algebraically closed field of
characteristic $p$.
Let $G$ be a linear algebraic group which is connected and
reductive. 
We fix  a maximal torus $T$ of $G$ of dimension $n$, 
the rank of $G$.
We also fix $B$, a Borel subgroup of $G$ with $B \supseteq T$ 
and let $W$ be the Weyl group of $G$.

We will write $\Mod(G)$ for the category of finite
dimensional rational $G$-modules. Most $G$-modules considered in this
paper will belong to this category. 
Let $X(T)=X$ be the weight lattice for $G$ 
 and $Y(T)=Y$ the dual weights.
The natural pairing $\langle -,- \rangle :X  \times Y
\rarr \Z$ is bilinear and induces an
isomorphism $Y\cong \Hom_\Z( X,\Z)$.
We take $R$ to be the roots of $G$. 
For each $\alpha \in R$ we take 
$\alphac \in Y$ to be  the coroot of $\alpha$. 
Let  $R^+$ be the positive roots, chosen so that $B$ is the negative
 Borel and let $S$ be the set of simple roots. 
Set $\rho = \half \sum_{\alpha \in R^+} \alpha \in
X\otimes_{\Z}\Q$.

We have a partial order on $X$ defined by 
$\mu \le \lambda \lrimpl \lambda -\mu \in \N S$.
A weight $\lambda$ is \emph{dominant} if
$\langle \lambda, \alphac \rangle \ge 0$ for all $\alpha \in S$ and
we let $X^+$ be the set of dominant weights.

Take $\lambda \in X^+$ and let $k_\lambda$ be the one-dimensional module
for $B$ which has weight $\lambda$. We define the induced
module, $\nabla(\lambda)= \Ind_B^G(k_\lambda)$. 
This module has formal character given by Weyl's character formula and has
simple socle $L(\lambda)$, 
the irreducible $G$-module of highest weight
$\lambda$.  Any finite dimensional, rational irreducible $G$-module
is isomorphic to $L(\lambda)$ for a unique $\lambda \in X^+$.

Since $G$ is split, connected and reductive we have an
antiautomorphism, $\tau$, which acts as the identity on $T$
(\cite{jantz}, II, corollary 1.16).
From this morphism we may define $^\circ$, 
a contravariant dual. It does not change a module's character, hence
it fixes the irreducible modules. 
We define the Weyl module, to be 
$\Delta(\lambda)=\nabla(\lambda)^\circ$.
Thus $\Delta(\lambda)$ has simple head $L(\lambda)$.

We return to considering the weight lattice $X$ for $G$.
There are also the affine reflections
$s_{\alpha,mp}$ for $\alpha$ a positive
root and $m\in \Z$ which act on $X$ as
$s_{\alpha,mp}(\lambda)=\lambda -(\langle\lambda,\alphac\rangle -mp )\alpha$.
These generate the affine Weyl group $W_p$. 
We mostly use the dot action of $W_p$ on $X$ which is 
the usual action of $W_p$, with the origin
shifted to $-\rho$. So we have $w \cdot \lambda = w(\lambda+\rho)-\rho$.
%
%
If $F$ is an alcove for $W_p$ then its closure $\bar{F}\cap X$ is a
fundamental domain for $W_p$ operating on $X$. The group $W_p$
permutes the alcoves simply transitively.
We set 
$C= \{ \lambda \in X \otimes _{\Z} \R \ \mid\  0< \langle \lambda +\rho,
\alphac \rangle < p \quad \forall\, \alpha \in R^+ \}$ 
 and call $C$ the \emph{fundamental alcove}.
We also set
$h = \max \{ \langle \rho, \betac \rangle +1 \ \mid\ \beta \in
R^+\}$.
When $R$ is irreducible then $h$ is the Coxeter number of $R$.
In general, it is the maximum of all Coxter numbers of the irreducible
components of $R$. 
We have
$C \cap X \ne \emptyset \ \lrimpl \ \langle \rho, \betac \rangle < p
\quad \forall\, \beta \in R^+ \ \lrimpl \ p \ge h$.

A facet $F$ is a \emph{wall} if 
there exists a unique $\beta\in R^+$ with $\langle \lambda
+\rho, \betac \rangle =mp$ for some $m \in \Z$ and for all $\lambda \in F$.
Let $s_F= s_{\beta, mp}$. 
This is the unique reflection in $W_p$
which acts as the identity on $F$ and we call $s_F$ the
reflection with respect to $F$.

Let 
$\Stab_{W_p}(\lambda)$ be all the elements of $W_p$ which stabilise
$\lambda\in X$. 
We take $\Sigma$ to be the set of all reflections $s_F$ where $F$
is a wall (for $W_p$) with $F \subset \bar{C}$. 
Thus the set $\Sigma$ consists of the reflections $s_{\alpha ,0}$ with
$\alpha \in S$ together with $s_{\beta,p}$ with $\beta$ the longest
short root of each irreducible component of the root system $R$.
Let $\Sigma^0(\mu)$ be
the subset of $\Sigma$ where each element of $\Sigma^0(\mu)$ fixes
$\mu$. 
The affine Weyl group $W_p$ is generated by $\Sigma$. These generators
form a presentation for $W_p$ as a Coxeter group so we may define a
length function $l(w)$ for $w \in W_p$ which is the 
length of a reduced expression
for $w$ in terms of elements of $\Sigma$.

We say that $\lambda$ and $\mu$ are \emph{linked} if  they belong to the 
same $W_p$ orbit on $X$ (under the dot action).
If two irreducible modules $L(\lambda)$ and $L(\mu)$ are in the same $G$
block then $\lambda$ and $\mu$ are linked. 

The category of rational $G$-modules has enough injectives and so we may
define $\Ext_G^*(-,-)$ as usual by using injective
resolutions (see \cite{benson}, section 2.4 and 2.5).
We will usually just write $\Ext$ for $\Ext_G$.

%

\section{Quasi-Hereditary Algebras}\label{sect:qha}
In this section we prove some lemmas about module category $\Mod(A)$,
for a quasi-hereditary algebra $A$ with poset $(\Lambda, \le)$,
standard modules $\Delta(\lambda)$ and costandard modules
$\nabla(\lambda)$. We will later lift these results to $\Mod(G)$. 

We say $X \in \Mod(A)$ has a \emph{good filtration} if it
has a filtration
$0=X_0 \subset X_1 \subset \cdots \subset X_i=X$ 
with quotients $X_j/X_{j-1}$ isomorphic to 
$\nabla(\mu_j)$ for some $\mu_j \in \Lambda$. The class of
$A$-modules with good filtration is denoted $\good$, and dually the class of
modules filtered by $\Delta(\mu)$'s is denoted $\doog$. We say that $X
\in \doog$ has a \emph{Weyl filtration}. The multiplicity of
$\nabla(\mu)$ in a filtration of $X \in \good$ is independent of the
filtration chosen and  is denoted by $\bigl(X:\nabla(\mu)\bigr)$. The composition
multiplicity of $L(\mu)$ in $ X\in \Mod(S)$ is denoted by $\big[X:L(\mu)\big]$.
We point out that even when $\nabla(\mu)=L(\mu)$ then it is still not
necessarily true that $\bigl(X:\nabla(\mu)\bigr)=\big[X:L(\mu)\big]$.

Some of the important properties of $\good$ and $\doog$ are stated
below.

\begin{propn}
\noindent
\begin{enumerate}
\item[(i)]\ 
Let $X \in \Mod(A)$ and $\lambda \in \Lambda$. If $\Ext_A^1\bigl(X,
\nabla(\lambda)\bigr)\ne 0$ then $X$ has a composition factor $L(\mu)$ with
$\mu > \lambda$.
\item[(ii)]\label{propn:good}\
For $X \in \doog$, $Y\in \good$ and $i>0$, we have 
$\Ext_A^i(X,Y)=0$.
\item[(iii)]\
Suppose $\Ext_A^1\bigl(\Delta(\mu),M\bigr)=0$ for all $\mu \in \Lambda$
then $M\in \good$.
\item[(iv)]\
Let $X \in \good$ (resp. $X \in \doog $) and $Y$ a direct summand of
$X$ then $Y\in \good$ (resp. $Y \in \doog$).
\end{enumerate}
\end{propn}
\begin{proof}
See
~\cite{donkbk}, A2.2.
\end{proof}


Suppose $X \in \Mod (A)$. We can resolve $X$ by modules $M_i\in \good$ as
follows
$$0\rarr X \rarr M_0 \rarr M_1 \rarr \dots \rarr M_d \rarr 0.$$
Such a resolution a \emph{good resolution} for $X$.
Good resolutions exist for all $A$-modules as $A$ has enough
injectives and an injective
resolution is also a good resolution.

The following definition may be found in~\cite{fripar} where 
a proof of the equivalence of properties (i) and (ii) may be found 
(see~\cite{fripar}, proposition 3.4).

\begin{defn}\label{defn:gfd}
Let $X \in \Mod (A)$. We say $X$ has \emph{good filtration dimension}
$d$, denoted  $\gfd(X)=d$,
if the following two equivalent conditions hold:
\begin{enumerate}
\item[(i)]
$0\rarr X \rarr M_0 \rarr M_1 \rarr \dots \rarr M_d \rarr 0$
is a resolution for $X$ with $M_i \in \good$, of shortest
possible length.
\item[(ii)]
$\Ext_A^i(\Delta(\lambda),X)=0$ for all $i>d$ and all $\lambda \in
\Lambda$,
but there exists $\lambda \in \Lambda$ such that
$\Ext_A^d(\Delta(\lambda),X)\ne 0$.
\end{enumerate}
\end{defn}

Similarly we have the dual notion of the 
\emph{Weyl filtration dimension} of $M$ which we will denote 
$\wfd(M)$. 

\begin{lem}\label{lem:wandg}
Given $A$-modules $M$ and $N$, we have
$$\Ext_A^i(N,M)=0\mbox{ for }i > \wfd(N)+\gfd(M).$$
\end{lem}
\begin{proof} 
See~\cite{parker1}, lemma 2.2.
\end{proof}

\begin{defn}
Let $g=\sup\{ \gfd(X) \mid X \in \Mod(A)\}$. 
We say $A$ has good filtration
dimension $g$ and denote this by $\gfd(A)=g$.
Let $w=\sup\{ \wfd(X) \mid X \in \Mod(A)\}$. 
We say $A$ has Weyl filtration
dimension $w$ and denote this by $\wfd(A)=w$.
\end{defn}

\begin{rem} In general $\gfd(A)$ is not the good filtration dimension
of $A$ when considered as its own left (or right) module. Similar
remarks apply to $\wfd(A)$. We will only
use $\gfd(A)$ and $\wfd(A)$ in the sense that they are defined above.
\end{rem}

For a finite dimensional $k$-algebra $A$, 
the \emph{injective dimension} of an $A$-module $M$, 
is the length of a shortest possible injective
resolution and is denoted by $\inj(M)$.
Equivalently we have
$\inj(M)=\sup\{ d \mid \Ext_A^d(N,M)\not\cong 0 \mbox{ for } N \in \Mod(A)\}.$
 The \emph{global dimension} of $A$ is the
supremum of all the injective dimensions for $A$-modules, and is
denoted by $\glob(A)$. This is equivalent to
$\glob(A)=\sup\{ d \mid \Ext_A^d(N,M) \not\cong 0 \mbox{ for some } N,M \in \Mod(A)\}.$
We will also denote the \emph{projective dimension} of an $A$-module
$M$ by $\proj(M)$.

\begin{cor}\label{cor:glob}
The global dimension of
$A$ has an upper bound of $\wfd(A)+\gfd(A)$.
\end{cor}

\begin{defn}
We say a module $T$ is a \emph{tilting module} if $T$ has both a good
filtration and a Weyl filtration. That is $T \in \good \cap
\doog$.
\end{defn}
For each $\lambda \in \Lambda$ there is a unique indecomposable
tilting module, $T(\lambda)$, of highest weight $\lambda$ with
$[T(\lambda):L(\lambda)]=1$. Every tilting module $T$ can be written
as a direct sum of indecomposable tilting modules $T(\mu)$ with $\mu
\in \Lambda$~(\cite{donkbk}, theorem A4.2).

\begin{defn}
Take $\lambda\in \Lambda$. 
We take a chain $\mu_0 < \mu_1 < \cdots < \mu_{l-1}< \mu_l=\lambda$ 
with $l$ maximal and $\mu_i \in \Lambda$.
We define the length of $\lambda$, $l(\lambda)$ to be $l$.
We also define $l(\Lambda)=\max\{l(\lambda)\mid \lambda \in
\Lambda\}$.
\end{defn}

\begin{lem}\label{lem:wfdlen}
We have $\wfd\bigl(\nabla(\lambda)\bigr) \le l(\lambda)$.
\end{lem}
\begin{proof}
If $l(\lambda)=0$ then $\lambda$ is minimal so
$\nabla(\lambda)=\Delta(\lambda)$
and $\wfd\bigl(\nabla(\lambda)\bigr)=0=l(\lambda)$.

Now suppose the lemma is true for $\mu < \lambda$.
We have a \ses
$$ 0\rarr N \rarr T(\lambda)\rarr \nabla(\lambda) \rarr 0$$
where $T(\lambda)$ is the indecomposable tilting module of highest weight
$\lambda$.
Applying $\Ext_A^*\bigl(-,\nabla(\nu)\bigr)$ for $\nu \in \Lambda$ gives us
$$ \Ext_A^{i-1}\bigl(N,\nabla(\nu)\bigl) \rarr 
\Ext_A^i\bigl(\nabla(\lambda),\nabla(\nu)\bigr)
\rarr \Ext_A^i\bigl(T(\lambda),\nabla(\nu)\bigr) $$
We now take $i > l(\lambda)$. All the $\nabla(\mu)$ appearing in a
good filtration of $N$  have $\mu < \lambda$. Hence $l(\mu) <
l(\lambda)$
and so $i-1> l(\mu)$. We now use the induction hypothesis to get
$\Ext_A^{i-1}\bigl(N, \nabla(\nu)\bigr) =0$.  We also have
$\Ext_A^i\bigl(T(\lambda),\nabla(\mu)\bigr)=0$ as $T(\lambda)$ is
tilting and using proposition~\ref{propn:good} (ii).
Hence $\Ext_A^i\bigl(\nabla(\lambda),\nabla(\mu)\bigr)=0$ for $i > l(\lambda)$ 
and $\wfd\bigl(\nabla(\lambda)\bigr)\le l(\lambda)$.
\end{proof}

We may similarly prove that $\gfd\bigl(\Delta(\lambda)\bigr)\le l(\lambda)$.

\begin{lem}\label{lem:wfdS}
$$\wfd{A}=\max\{\wfd\bigl(\nabla(\lambda)\bigr)\mid \lambda\in \Lambda\}.$$
\end{lem}
\begin{proof}
We certainly have
$$\wfd{A}\ge\max\{\wfd\bigl(\nabla(\lambda)\bigr)\mid \lambda\in \Lambda\}.$$
Take $\lambda\in \Lambda$ with $\wfd\bigl(L(\lambda)\bigr)=d=\wfd(A)$.
Let $Q$ be the quotient $\nabla(\lambda)/L(\lambda)$.
Since $\wfd\bigl(L(\lambda)\bigr)$ was maximal we must have $\wfd(Q) \le d$.
Let $\mu\in\Lambda$, the corresponding \les for the \ses for $Q$ gives us
$$
\Ext_A^d\bigl(\nabla(\lambda),\nabla(\mu)\bigr) \rarr
\Ext_A^d\bigl(L(\lambda),\nabla(\mu)\bigr) \rarr
\Ext_A^{d+1}\bigl(Q,\nabla(\mu)\bigr)
$$
Now $\Ext_A^{d+1}\bigl(Q,\nabla(\mu)\bigr)=0$ for all $\mu\in\Lambda$ 
by lemma~\ref{lem:wandg}.
But there exists $\mu\in \Lambda$ with
$\Ext_A^d\bigl(L(\lambda),\break\nabla(\mu)\bigr)\ne 0$.
Hence $\Ext_A^d\bigl( \nabla(\lambda),\nabla(\mu)\bigr)\ne 0$.
Thus there exists $\lambda\in\Lambda$ with $\wfd\bigl(\nabla(\lambda)\bigr)=d
=\wfd(A)$
\end{proof}

\begin{rem}\label{rem:wfdS}
We may replace the set of $\nabla(\lambda)$ with any set of
$A$-modules $\mathcal{X}$ with the property that for all $\lambda \in
\Lambda$ there exists $X \in
\mathcal{X}$ with $L(\lambda)$ contained in the socle of $X$.
We can then repeat the argument above to get
$\wfd(A)=\max\{\wfd(X)\mid X \in \mathcal{X}\}$.
\end{rem}

Now suppose $A$ is a quasi-hereditary algebra with contravariant
duality preserving simples.
That is there exists an involutory, contravariant functor 
$^\circ: \Mod (A)\rarr \Mod (A)$
such that, $\Delta(\lambda)^\circ \cong \nabla(\lambda)$ 
(and $\Ext_A^i(M,N) \cong \Ext_A^i(N^\circ ,M^\circ)$). 
We will usually shorten this and say $A$ has a simple preserving
duality.

\begin{rem}It is clear (given the equivalences in the definition for the
good filtration dimension) that for $A$ with simple preserving duality
and $M$ an $A$-module 
we have $\wfd(M) = \gfd(M^\circ)$.
We will use this without further comment.
\end{rem}

Thus lemma~\ref{lem:wfdlen} gives an upper bound
for $\wfd(A)$ of $l(\Lambda)$.
Corollary~\ref{cor:glob} gives, for $A$ with simple preserving duality
that
$$\glob(A)\le 2\gfd(A)=2\wfd(A)\break \le 2l(\Lambda).$$

We say a subset $\Pi$ of a poset $(\Lambda, \le)$ (not necessarily finite) is
\emph{saturated} if for all $\lambda \in \Pi$ then $\mu \le \lambda$
implies that $\mu \in \Pi$. 

Take $G$ to be a split, connected reductive algebraic group with
weight lattice $X$.
Suppose $\Pi$ is a finite saturated subset of $X^+$ with respect to
the dominance ordering. We may consider
$G$-modules all whose composition factors have highest weights lying
in $\Pi$. These modules form a subcategory of $\Mod(G)$  
which is a highest weight category corresponding to a quasi-hereditary
algebra which we denote $S(\Pi)$, the generalised Schur algebra
(see~\cite{donk1} for more
information). We have a natural isomorphism
$$\Ext_{S(\Pi)}^i(M,N) \cong  \Ext_G^i(M,N)$$
for $S(\Pi)$-modules $M$ and $N$~\cite{donk1}, 2.2d. 
The costandard and standard modules for $S(\Pi)$ are exactly the
induced and Weyl modules for $G$ respectively.
Thus as long as we restrict our attention to finite dimensional
$G$-modules then we can lift the results from quasi-hereditary
algebras to $G$.

Generally speaking, a finite dimensional $G$-module does not have a
finite injective or projective resolution. It will have, however, have a finite good
(and Weyl) resolution. Thus we can lift the definitions of good (and Weyl)
filtration dimension to $\Mod(G)$. 

If we take $G=\GL_n(k)$ and $\Pi=\partn(n,r)$ then $S(\Pi)$ is
isomorphic to $S(n,r)$, the usual Schur algebra. Thus
Schur algebras are quasi-hereditary with poset $\partn(n,r)$ ordered
by dominance. 


\section{Properties of Translation Functors}\label{sect:trans}
For any $G$-module $V$ and any $\mu\in X$, set
$\pr_\mu V$ equal to the sum of submodules of $V$ such that all the
composition factors have highest weight in $W_p \cdot \mu$. Then
$\pr_\mu V$ is the largest submodule of $V$ with this property.
The following definition is due to Jantzen~\cite{jantz}, II, 7.6.

\begin{defn}
Suppose $ \lambda$, $\mu\in \bar{C}$. There is a unique $\nu_1 \in
X^+\cap W(\mu-\lambda)$. We define the \emph{translation functor}
$T_\lambda ^\mu$ from $\lambda$ to $\mu$ via
$$T_\lambda ^\mu V = \pr _\mu (L(\nu_1)\otimes \pr_\lambda V)$$
for any $G$-module $V$. It is a functor from $\Mod(G)$ to
itself.
\end{defn}

\begin{lem}\label{lem:adjoint}
Let $\lambda$ and $\mu \in \bar{C}$,
then the functors $T_{\lambda}^\mu$ and $T_{\mu}^\lambda$ are adjoint to
each other. 
For $M,N \in \Mod(G)$ we have
$\Ext^i(T_\lambda^\mu M,N) \cong \Ext^i(M,T_\mu^\lambda N)$.
\end{lem}
\begin{proof}
See~\cite{jantz}, II, lemma 7.6 (b) and remark 7.6 (2).
\end{proof}

\begin{propn}\label{propn:transgood}
Let $\mu$, $\lambda\in \bar{C}$ and $w \in W_p$ with $w\cdot \mu
\in X^+$, then $T_\mu ^\lambda \nabla(w \cdot \mu) $ has a good
filtration.
Moreover the factors are $\nabla(ww_1 \cdot \lambda)$ with $w_1 \in
\Stab_{W_p}(\mu)$ and $ww_1\cdot\lambda \in X^+$. Each different $ww_1
\cdot \lambda$ occurs exactly once.
\end{propn}
\begin{proof}
See~\cite{jantz}, proposition 7.13.
\end{proof}

\begin{cor}\label{cor:transses}
Let $\lambda\in C$ and $\mu \in \bar{C}$. Suppose there is
$s \in \Sigma$ with $\Sigma^0(\mu)=\{s\}$.
Let $w \in W_p$ with
$w\cdot \lambda\in X^+$ and $w\cdot \lambda < ws\cdot \lambda$.
Then we have a \ses 
$$ 0\rarr \nabla(w\cdot \lambda)
\rarr T_\mu^\lambda\nabla(w\cdot \mu)
\rarr \nabla(ws\cdot \lambda)\rarr 0.$$
\end{cor}
\begin{proof}
See~\cite{jantz}, lemma 7.19 (a).
\end{proof}

We would like to know when such a situation in the above corollary
occurs.
Firstly we need a $\lambda \in C$ and  this happens when $p \ge h$.
We also need a weight $\mu$ lying
on the wall between $\lambda$ and $s\cdot \lambda$. This happens when
the derived group of $G$ is simply connected and $p \ge h$.
See~\cite{jantz}, II, 6.3 (1), for details.
We will henceforth assume that $p \ge h$ and that the derived group of
$G$ is simply connected.
We will also assume  that the root system $R$ of $G$
is irreducible, although we believe that theorem~\ref{thm:wfdlamb} 
is also true in the more general case.


We have another  partial order on $X$ denoted $\up$.
If $\alpha$ is a positive root  and $m\in \Z$ then we set
$$s_{\alpha,mp} \cdot \lambda \up \lambda \quad\mbox{if and only
if}\quad \langle\lambda+\rho ,
\alphac\rangle \ge mp.$$
This then generates an order relation on $X$.
So $\mu \up \lambda$ if there are reflections $s_i \in W_p$ with 
$$\mu=s_m s_{m-1} \cdots s_1 \cdot \lambda \up 
s_{m-1} \cdots s_1 \cdot \lambda \up 
\cdots \up
s_1 \cdot \lambda \up 
\lambda .$$

We define $l(\lambda)$ for $\lambda\in X^+$ to be the length of a maximal chain 
$\mu_0 \up \mu_1 \up \cdots \up \mu_{l-1} \up \mu_l=\lambda$
with $\mu_0 \in \bar{C}$, each $\mu_i \ne \mu_{i+1}$ and $\mu_i \in X$.
We will also define $\bl(\lambda)$ for $\lambda \in X^+$ to be the
length of a maximal chain 
$\mu_0 \up \mu_1 \up \cdots \up \mu_{l-1} \up \mu_l=\lambda$
with all $\mu_i \in X^+$. 

We define $d(\lambda)$ 
to be the number of hyperplanes separating $\lambda$ and  a weight
lying in $C$ (we do not count any hyperplanes that $\lambda$ may lie
on).
Take $n_\alpha$, $d_\alpha \in \Z$ with $\langle \lambda+\rho,
\alphac \rangle = n_\alpha p +d_\alpha$ and $0< d_\alpha \le p$
for all $\alpha$ a positive root. 
If $\lambda$ is dominant then $d(\lambda) = \sum_{\alpha>0} n_\alpha$. 

\begin{lem}\label{lem:alclen}
If $\lambda \in C$ and $w \in W_p$ with $\wlamb \in X^+$ 
then $\bl(\wlamb)=l(\wlamb)= l(w)=d(\wlamb)$
\end{lem}
\begin{proof}
Since $\wlamb$ lies inside an alcove we have that 
$d(\wlamb)=l(w)$. (This is true as the alcoves in $X$ can be identified
with chambers in the Coxeter complex associated to $W_p$.) 
It is clear that $l(\wlamb) \ge \bl(\wlamb)$.
We have using~\cite{jantz}, proposition 6.8, that 
$\bl(\wlamb) \ge d(w\cdot \lambda)$.
Now take a maximal chain for $\wlamb$,
$\mu_0 \up \mu_1 \up \cdots \up \mu_l=\wlamb$ with $\mu_0 \in C$ and
$\mu_i \in X$.
We know that in this chain for $\wlamb$ 
we have $d(\mu_i) < d(\mu_{i+1})$, by applying~\cite{jantz}, lemma 6.6.
Thus $d(\wlamb) \ge l(\wlamb)$. 
Hence we have the equalities as claimed.
\end{proof}

\begin{rem}\label{rem:bruhat}
If $\lambda \in C$ then the $\up$-ordering on $X^+\cap W_p \cdot
\lambda$ is equivalent to the Bruhat ordering on $W_p$.
That is we have for $\lambda \in C$ and $w ,v \in W_p$ with $w \cdot
\lambda$ and $v \cdot \lambda \in X^+$ that 
$$ w \cdot \lambda \up v \cdot \lambda\quad\mbox{if and only if}\quad
w \le v .$$
This can be seen from the definition of the Bruhat order in
\cite{humph2}, section 5.9, and using the previous lemma.
See also~\cite{verma}, section 1.6.
\end{rem}

We have that $[\nabla(\lambda):L(\mu)] \ne 0$ implies $\mu \up
\lambda$~\cite{ander1}, corollary 3, (known as the strong linkage principle).
Thus when we take 
$\Pi$, a finite saturated subset of $X^+$ with respect to
the $\up$ ordering, the corresponding algebra $S(\Pi)$ is
quasi-hereditary, 
thus we may apply lemma~\ref{lem:wfdlen} to deduce that
$\wfd\bigl(\nabla(\lambda)\bigr)\le \bl(\lambda)$.


\section{The Weyl Filtration Dimension of the Induced
Modules}\label{sect:wfd}
\begin{lem}\label{lem:wfdmu}
Suppose we have the situation of corollary~\ref{cor:transses}.
So we have $\lambda\in C$, $\mu \in \bar{C}\backslash C$, $w\cdot \lambda<
ws\cdot \lambda$, $w\cdot \lambda \in X^+$ and $\Sigma^0(\mu)=\{s\}$.
If $l(w) \ge 1$ then $\wfd\bigl(\nabla(w \cdot \mu)\bigr)< l(w)$.
\end{lem}
\begin{proof}
It is clear that any non-repeating chain for $w\cdot \mu$,
$w_1\cdot\mu \up \cdots \up w_i\cdot \mu \up \cdots  \up w_m\cdot
\mu=w\cdot\mu$
with $w_i\cdot \mu \in X^+$ gives a non-repeating  chain 
$w_1\cdot \lambda \up \cdots\up  w_i \cdot \lambda\up 
\cdots \up w_m \cdot \lambda=w\cdot \lambda$ with $w_i\cdot \lambda \in X^+$.
So $\bl(\wmu) \le \bl(\wlamb)=l(\wlamb)$.

If $m=l(w)$ then the chain for $\lambda$ is maximal. So we would have
$w_1= 1$ and $w_2= s_{\beta,p}$ for $\beta$ the longest short root of
$R$ (as $R$ is irreducible).
But then $w_1\cdot \mu = \mu \in X^+$. We also assumed $\mu \in
\bar{C} \backslash C$ so $\mu$ must be fixed by $s_{\beta,p}$.
So we have $\mu = w_1 \cdot \mu = w_2 \cdot \mu$. But this means
the chain for $\mu$ repeats -- a contradiction.
Thus $\bl(\wmu) < l(\wlamb)=l(w)$ by lemma~\ref{lem:alclen}.
Now lemma~\ref{lem:wfdlen} gives us the result.
%
%
\end{proof}

\begin{thm}\label{thm:wfdlamb}
Suppose the root system $R$ of $G$ is irreducible 
and $\lambda \in C$. Then
$$\wfd\bigl(\nabla(\wlamb)\bigr)= l(w).$$ 
\end{thm}
\begin{proof}
We proceed by induction on $l(w)$.
If $l(w)=0$ then $\nabla(\lambda)=\Delta(\lambda)=L(\lambda)$ so
$\wfd\bigl(\nabla(\lambda)\bigr)=0$. 

Now let $w=s$, $s \in \Sigma$ with $s\cdot \lambda \in X^+$.
Take $\mu$ to be a dominant weight on the wall separating 
$\lambda$ and $s\cdot \lambda$. Such a $\mu$ has
$\wfd\bigl(\nabla(\mu)\bigr)=0$. 
Thus $T_\mu^\lambda\bigl(\nabla(\mu)\bigr)$ is a tilting
module of highest weight $s \cdot \lambda$.
So the \ses of corollary~\ref{cor:transses} is a Weyl resolution of
$\nabla(s \cdot \lambda)$ and so $\wfd\bigl(\nabla(s\cdot \lambda)\bigr)\le 1$.
But $\Ext^1\bigl(\nabla(s \cdot \lambda), \nabla(\lambda)\bigr)=k$ 
by~\cite{jantz}, II, proposition 7.21, and so $\wfd\bigl(\nabla(s\cdot
\lambda)\bigr)= 1=l(s)$. 

Now suppose the theorem is true for all $w \in W_p$ with $l(w)\le
l$, $l\ge 1$. We will show the result holds for $ws$ with $s \in
\Sigma$. We take $\mu \in \bar{C}\backslash C$ with
$\Sigma_0(\mu)=\{s\}$.
We have for all $i$, $v \in W_p$ and $v \cdot \lambda \in X^+$
$$\Ext^i\bigl(T_\mu^\lambda\bigl(\nabla(\wmu)\bigr), \nabla(v \cdot
\lambda )\bigr)
\cong \Ext^i\bigl(\nabla(\wmu), T_\lambda^\mu(\nabla(v \cdot \lambda)\bigr)\bigr)
\cong \Ext^i\bigl(\nabla(\wmu), \nabla(v \cdot \mu)\bigr)$$
by lemma~\ref{lem:adjoint} and proposition~\ref{propn:transgood}.
So we have 
\begin{equation}\label{wfdlen}
\wfd\bigl(T_\mu^\lambda\bigl(\nabla(\wmu)\bigr)\bigr)= 
\wfd\bigl(\nabla(\wmu)\bigr) < l(w)
\end{equation}
by lemma~\ref{lem:wfdmu}.

Applying $\Ext^*\bigl(-,\nabla(\nu)\bigr)$ with $\nu \in X^+$ to the \ses of
corollary~\ref{cor:transses} gives us
\begin{multline*}
\Ext^{i}\bigl(T_\mu^\lambda\bigl(\nabla(\wmu)\bigr),\nabla(\nu)\bigr)
\rarr
\Ext^i\bigl(\nabla(\wlamb),\nabla(\nu)\bigr) 
\\ 
\rarr
\Ext^{i+1}\bigl(\nabla(\wslamb),\nabla(\nu)\bigr) 
\rarr
\Ext^{i+1}\bigl(T_\mu^\lambda\bigl(\nabla(\wmu)\bigr),\nabla(\nu)\bigr).
\end{multline*}
Thus for $i \ge l(w)$ we have
$$
\Ext^i\bigl(\nabla(\wlamb),\nabla(\nu)\bigr) 
\cong
\Ext^{i+1}\bigl(\nabla(\wslamb),\nabla(\nu)\bigr) 
$$
using \eqref{wfdlen} and lemma~\ref{lem:wandg}.
Hence $\wfd\bigl(\nabla(\wslamb)\bigr)= \wfd\bigl(\nabla(\wlamb)\bigr)+1=l(w)+1=l(ws)$,
as required.
\end{proof}

We may use the $^\circ$-duality to get that
$\gfd(\Delta(\wlamb))=l(w)$.
The previous theorem and lemma~\ref{lem:wandg} give us that 
for $v \in W_p$ with $\vlamb \in X^+$ we have
$\Ext^{i}\bigl(\nabla(\wlamb),\Delta(\vlamb)\bigr)=0$ for $i >
l(w)+l(v)$.
The following corollary tells us that this bound is strict.

\begin{thm}\label{thm:extiind}
Suppose $\lambda \in C$, and $w$, $v \in W_p$ with $\wlamb$, $\vlamb \in X^+$. Then
$$\Ext^{l(w)+l(v)}\bigl(\nabla(\wlamb),\Delta(\vlamb)\bigr)\cong k.$$
\end{thm}
\begin{proof}
We proceed by induction on $l(w)+l(v)$.
If $l(w)+l(v)=0$ then $\wlamb=\vlamb=\lambda$ so
$\Hom\bigl(\nabla(\wlamb),\Delta(\vlamb)\bigr)\cong
\Hom\bigl(\nabla(\lambda),\Delta(\lambda)\bigr)\cong k$.

If $l(w)+l(v)=1$ then either $\wlamb= \lambda$ or $\vlamb = \lambda$.
Without loss of generality (using the $^\circ$-duality), take
$\wlamb \ne \lambda$. Thus $\Delta(\vlamb)=\nabla(\lambda)$. Also 
$l(w)=1$ so $w=s \in \Sigma$.
By~\cite{jantz}, II, proposition 7.21 (c), we have
$\Ext^1\bigl(\nabla(s \cdot \lambda), \nabla(\lambda) \bigr) \cong k$.
Thus the corollary is true for $l(w)+l(v)=1$.

Now take $l(w)=l(v)=1$ 
Applying $\Ext^*\bigl(\nabla(s \cdot \lambda),-\bigr)$ to 
the $^\circ$-dual of the \ses of
corollary~\ref{cor:transses} gives us
\begin{equation*}
\Ext^1\bigl(\nabla(s \cdot \lambda),T_\mu^{\lambda}\bigl(\Delta(\mu)\bigr)\bigr)
\rarr
\Ext^1\bigl(\nabla(s \cdot \lambda),\Delta(\lambda)\bigr)
\rarr
\Ext^2\bigl(\nabla(s \cdot \lambda),\Delta(s \cdot \lambda)\bigr)
\rarr 0.
\end{equation*}
The last zero follows by lemma~\ref{lem:wandg}. Also 
$$\Ext^1\bigl(\nabla(s \cdot
\lambda),T_\mu^{\lambda}\bigl(\Delta(\mu)\bigr)\bigl)
\cong 
\Ext^1\bigl(T^\mu_{\lambda}(\nabla(s \cdot \lambda)),\nabla(\mu)\bigr)
\cong 
\Ext^1\bigl(\nabla(\mu),\nabla(\mu)\bigr)
\cong
0.
$$
Hence 
$$
\Ext^2\bigl(\nabla(s \cdot \lambda),\Delta(s \cdot \lambda)\bigr)
\cong
\Ext^1\bigl(\nabla(s \cdot \lambda),\Delta(\lambda)\bigr)
\cong k.
$$

Now suppose the corollary is true for all $w, v \in W_p$ with $\wlamb,
\vlamb \in X^+$ and $l(w)+l(v)\le m$, for some $m\ge 1$.
We need to show the result holds for $l(w^\prime)+l(v^\prime)=m+1$,
$w^\prime, v^\prime \in W_p$ and $w^\prime \cdot \lambda, v^\prime
\cdot \lambda \in X^+$.
Without loss of generality we may take 
$v^\prime = v$ and $w^\prime =ws$ with $s \in \Sigma$.
We may also assume that $l(v^\prime)$ or $l(w^\prime)$ is at least 2
so that we can assume $w \ne 1$. (As we have already covered the case
with $l(w)=l(v)=1$.)

Apply $\Ext^*\bigl(-,\Delta(\vlamb)\bigr)$ to the
\ses of corollary~\ref{cor:transses} to get 
\begin{multline*}
\Ext^{m}\bigl(T_\mu^{\lambda}\bigl(\nabla(\wmu)\bigr),\Delta(\vlamb)\bigr)
\rarr
\Ext^{m}\bigl(\nabla(\wlamb),\Delta(\vlamb)\bigr)
\\ 
\rarr
\Ext^{m+1}\bigl(\nabla(\wslamb),\Delta(\vlamb)\bigr)
\rarr
\Ext^{m+1}\bigl(T_\mu^{\lambda}\bigl(\nabla(\wmu)\bigr),\Delta(\vlamb)\bigr).
\end{multline*}
But $\wfd\bigl(T_\mu^\lambda\bigl(\nabla(\wmu)\bigr)\bigr)< l(w)$ by
\eqref{wfdlen} (provided $w \ne 1$).
Now we may apply lemma~\ref{lem:wandg} to get that the first and last
$\Ext$ groups above are zero. Thus the middle two groups are
isomorphic.
So by induction we have
\begin{equation*}
\Ext^{l(w^\prime) +l(v)}\bigl(\nabla(w^\prime \cdot \lambda),\Delta(\vlamb)\bigr)
\cong
\Ext^{l(w)+l(v)}\bigl(\nabla(\wlamb),\Delta(\vlamb)\bigr)
\cong k.\qedhere
\end{equation*}
\end{proof}

\begin{cor}\label{cor:wfdnabd}
For $\lambda, \mu \in X^+$ lying inside an alcove and in the same
$W_p$-orbit we have
$$\wfd\bigl(\nabla(\lambda)\bigr)= d(\lambda)\mbox{,}\quad
\gfd\bigl(\Delta(\mu)\bigr)= d(\mu)
\quad\mbox{and}\quad
\Ext^{d(\lambda)+d(\mu)}\bigl(\nabla(\lambda),\Delta(\mu)\bigr) \cong k.$$
\end{cor}
\begin{proof}
We have that $\lambda= w \cdot \lambda_0$ and $\mu=v \cdot \lambda_0$ 
for some $\lambda_0 \in C$.
Lemma~\ref{lem:alclen}, theorem~\ref{thm:wfdlamb} and the previous corollary
then give us the result.
\end{proof}

\begin{cor}\label{cor:wfdirr}
For $\lambda, \mu \in X^+$ lying inside an alcove and in the same
$W_p$-orbit we have
$$\wfd\bigl(L(\lambda)\bigr)= d(\lambda)
\quad\mbox{and}\quad
\Ext^{d(\lambda)+d(\mu)}\bigl(L(\lambda),L(\mu)\bigr) \cong k.$$
\end{cor}
\begin{proof}
Let $Q$ be the quotient $\nabla(\lambda)/L(\lambda)$.
If $L(\nu)$ is a composition factor of $Q$ then $\nu \up \lambda$ and
$\nu \ne \lambda$. Thus $l(\nu) < l(\lambda)$.
Hence $\wfd(Q) < l(\lambda)=d(\lambda)=l$. Now apply
$\Ext^*\bigl(-,\nabla(\nu)\bigr)$ to the \ses
$$0 \rarr L(\lambda) \rarr \nabla(\lambda) \rarr Q \rarr 0$$
to get
$$
\cdots \rarr \Ext^l\bigl(Q,\nabla(\nu)\bigr) \rarr
\Ext^l\bigl(\nabla(\lambda),\nabla(\nu)\bigr)
\rarr \Ext^l\bigl(L(\lambda),\nabla(\nu)\bigr)\rarr 0
$$
where the last zero follows by lemma~\ref{lem:wandg}.
We also have that $\Ext^l\bigl(Q,\nabla(\nu)\bigr)=0$ by
lemma~\ref{lem:wandg}.
Thus $\wfd\bigl(L(\lambda)\bigr)= d(\lambda)=l$ as required.

A similar argument yields that 
\begin{equation*}
\Ext^{d(\lambda)+d(\mu)}\bigl(L(\lambda),L(\mu)\bigr)
\cong \Ext^{d(\lambda)+d(\mu)}\bigl(\nabla(\lambda),\Delta(\mu)\bigr)\cong
k.\qedhere 
\end{equation*}
\end{proof}

The result of Ryom-Hansen's in the appendix, theorem 2.4, states that
for $\lambda \in C$ and  $w,v \in W_p$ with $v \le w$ and $\wlamb, 
\vlamb \in X^+$ 
$$ \Ext^{l(w)-l(v)}\bigl(L(w \cdot \lambda), \nabla(v\cdot
\lambda)\bigr) \cong k. $$
We also know that if $i >l(w) -l(v)$ then
$\Ext^i\bigl(L(w \cdot \lambda), \nabla(v\cdot\lambda)\bigr) \cong 0$
by the appendix, lemma 2.1, (see also \cite{jantz}, proposition 6.20).
So using this result and given remark~\ref{rem:bruhat} 
and lemma~\ref{lem:alclen} we may now prove

\begin{propn}\label{propn:nabnab}
Let $\lambda \in C$, $w,v \in W_p$ with $\wlamb, \vlamb \in X^+$
and $\vlamb \up \wlamb$ then 
$$\Ext^{l(w)-l(v)}\bigl(\nabla(\wlamb),\nabla(\vlamb)\bigr)
\cong k.$$ 
\end{propn}
\begin{proof}
We may argue along similar lines to the proof of
corollary~\ref{cor:wfdirr}.
\end{proof}

We now are in a position where we may deduce the projective and
injective dimensions of several modules for the generalised Schur
algebras.
We define $\Pi(\lambda)$ to be the (finite) saturated subset
of $X^+$ with respect to the $\up$-ordering whose highest weight is $\lambda$.

\begin{thm}\label{thm:injproj}
Suppose $\lambda\in X^+$ is regular (lies inside an alcove) then
in $\Mod(S(\Pi(\lambda)))$  for $\mu \in \Pi(\lambda)$ we have
$$
\inj(L(\mu))=\proj(L(\mu))=\proj(\nabla(\mu))=\inj(\Delta(\mu))=d(\mu)+d(\lambda)
$$
$$
\inj(\nabla(\mu))=\proj(\Delta(\mu))=d(\lambda)-d(\mu).
$$
\end{thm}
In particular this gives us information for the blocks of the Schur
algebra whose weights are regular.

\section{The Global Dimension of $S(n,r)$ when $p>n$}\label{sect:glob}
We will now focus on the classical Schur algebra. 
So $G=\GL_n(k)$, the root system of $\GL_n$ is
irreducible and its derived subgroup $\SL_n$ is simply connected.
We wish to determine the good filtration dimension and global
dimension for $S(n,r)$ (that is for the \emph{whole} Schur algebra,
not just for the regular blocks). 
So we need to know what $d(\lambda)$ is for 
$\lambda$ a partition and a condition for $\lambda$ to lie inside an
alcove. 
The next two lemmas do this.

\begin{lem}\label{lem:dlambA}
Suppose $G=\GL_n(k)$ and $\lambda =(\lambda_1,\lambda_2,\ldots,
\lambda_n)\in X^+$. Then we have
$$d(\lambda)= \sum_{i=1}^{n-1}\sum_{j=i+1}^{n}
               \left\lfloor \frac{\lambda_i -\lambda_j
                        -i+j-1}{p}\right\rfloor.$$
\end{lem}
\begin{proof}
Let $e_i=(0,\ldots,0,1,0,\ldots,0)\in X$ with a one in the $i$th
position. The $e_i$ form the usual basis of $X$, so $\lambda =
(\lambda_1,\lambda_2,\ldots , \lambda_n)=\sum_{i=1}^n \lambda_i e_i$.
We take $\omega_i=\sum_{j=1}^i e_j$. 
We can write $\lambda= \sum_{i=1}^{n-1} (\lambda_i-\lambda_{i+1})\omega_i
+\lambda_n\omega_n$.
Thus for $\alpha = e_i-e_j\in R^+$, we have
$$\langle \lambda+\rho, \alphac \rangle =\lambda_i-\lambda_j+j-i.$$
The definition of $d(\lambda)$ then gives us the result.
\end{proof}
\begin{lem}\label{lem:weightalc}
A weight $\lambda \in X^+$ lies inside an alcove if there exist
no integers $i$ and $j$ such that $\lambda_i -\lambda_j \equiv i-j
\pmod p$.
\end{lem}
\begin{proof}
A weight $\lambda \in X^+$ lies on a wall if there exists $\alpha \in
R^+$ such that $\langle \lambda +\rho, \alphac \rangle = mp$ for some
$m \in \Z$.
So a weight $\lambda$ does not lie on a wall if for all $\alpha=e_i-e_j$ we have 
$\lambda_i-\lambda_j +j-i \not\equiv 0 \pmod p$.
\end{proof}

We first calculate an upper bound for $\wfd\bigl(S(n,r)\bigr)$.

Let $E=L(1,0,\ldots,0)$ be the natural module for $\GL_n$.
We take $S^r E$ to be the $r$th symmetric power of $E$ and 
$\Wedge^r E$ to be the $r$th exterior power.
For $\lambda = (\lambda_1,\lambda_2, \ldots, \lambda_n)\in \partn(n,r)$ we take
$S^\lambda E= S^{\lambda_1} E \otimes S^{\lambda_2} E \otimes \cdots
\otimes S^{\lambda_n} E$ with $S^1 E = E$ and $S^0 E = k$.

\begin{lem}\label{lem:wfdS2}
$$\wfd\bigl(S(n,r)\bigr)=\max\{\wfd(S^\lambda E)\mid \lambda\in
\partn(n,r)\}.$$
\end{lem}
\begin{proof}
We have that $L(\lambda)$ embeds in $S^\lambda E$ 
by~\cite{donkbk}, section 2.1 (15)(i)(b).
So the set $\mathcal{X}=\{S^\lambda E \mid \lambda \in \partn(n,r)\}$
satisfies the requirements of remark~\ref{rem:wfdS}.
\end{proof}

For all $\lambda$ and $\mu \in X^+$ the module
$\nabla(\lambda) \otimes \nabla(\mu)$ has a good filtration.
A proof of this property, for type $A_n$, is given in~\cite{wang1}. 
It is proved for most other cases in~\cite{donkrat}. The general proof
is given in~\cite{mathieu}.
The $\nabla(\nu)$ which appear as quotients in this
filtration are given by Brauer's character formula~\cite{jantz}, II, lemma 5.8.
We can generalise this property to good and Weyl filtration dimensions
as below.

\begin{lem}\label{lem:wfdten}
Let $X$, $Y$ be $G$-modules
then we have
$$\wfd(X\otimes Y)\le \wfd(X)+ \wfd(Y).$$
\end{lem}
\begin{proof}
See~\cite{fripar}, proposition 3.4 (c), where the corresponding result
for good filtration dimensions is proved. 
\end{proof}

\begin{lem}\label{lem:ses}
We have the short exact sequence,
\begin{multline*}
0 \rarr \nabla(mp-j,1^{j},0^{n-j-1} ) \rarr 
S^{mp-j}E \otimes \Wedge^{j} E
\rarr \nabla(mp-j+1,1^{j-1},0^{n-j}) \rarr 0.
\end{multline*}
\end{lem}
\begin{proof}
Since $S^{mp-j}E \otimes \Wedge^{j} E$ has a good filtration
by the dual version of lemma~\ref{lem:wfdten}, this follows using characters.
\end{proof}
\begin{propn}\label{propn:wfdSless}
$$\wfd(S^rE) \le (n-1)\left\lfloor\frac{r}{p}\right\rfloor.$$
\end{propn}
\begin{proof}
We first reduce to the case $S^{mp}E$. Write
$r=r_0+pm$. If $0<r_0<p$ then the multiplication $S^{r_0}E\otimes S^{rm}E\to
S^rE$ splits~\cite{donkbk}, section 4.8, proposition (12), 
so that $\wfd(S^rE)\le
\wfd(S^{r_0}E\otimes S^{pm}E)\le \wfd(S^{r_0}E)+\wfd(S^{pm}E)=
\wfd(S^{pm}E)$ as $S^{r_0}(E) \in \doog$. 
So suppose $r=mp$.
We prove this proposition by induction on $m$.
The proposition is clearly true for $m=0$

We note that the modules $\Wedge^jE$ are tilting modules for
$S(n,r)$~\cite{donktilt}, lemma 3.4 (ii), and hence have Weyl 
filtration dimension $0$.
Dimension shifting using the induction hypothesis,
lemma~\ref{lem:wfdten},
lemma~\ref{lem:ses} and lemma~\ref{lem:wandg} gives us
\begin{align*}
\Ext^i\bigl(S^rE,\nabla(\mu)\bigr) &\cong
\Ext^{i-1}\bigl(\nabla(r-1,1,0^{n-2}),\nabla(\mu)\bigr) 
\\ 
&\cong \cdots\cong
\Ext^{i-j}\bigl(\nabla(r-j,1^j,0^{n-j-1}),\nabla(\mu)\bigr)
\end{align*}%
for $\mu\in \partn$ and $i-j>(n-1)(m-1)\ge\wfd(S^{mp-j}E\otimes
\Wedge^jE)$.
So for $i>(n-1)m$ we have
\begin{align*}
\Ext^i\bigl(S^rE,\nabla(\mu)\bigr)& \cong
\Ext^{i-n+1}\bigl(\nabla(r-n+1,1^{n-1}),\nabla(\mu)\bigr)
\\
&\cong\Ext^{i-n+1}\bigl(S^{r-n}E\otimes \Wedge^nE,\nabla(\mu)\bigr)
\cong 0
\end{align*}%
by induction and lemma~\ref{lem:wfdten}.
Hence $\wfd(S^rE) \le (n-1)m$ as required.
\end{proof}

Let $\lambda \in \partn$ with $\lambda=(\lambda_1,\lambda_2,\ldots,\lambda_n)$.
We define 
$$\left\lfloor\frac{\lambda}{p}\right\rfloor = 
\sum_{i=1}^n \left\lfloor\frac{\lambda_i}{p}\right\rfloor.$$

\begin{cor}
$$\wfd\bigl(S(n,r)\bigr) \le (n-1)\left\lfloor\frac{r}{p}\right\rfloor.$$
\end{cor}
\begin{proof}
We have $\wfd(S^\lambda E)\le \bigl\lfloor\frac{\lambda}{p}\bigr\rfloor$ 
using lemma~\ref{lem:wfdten} and proposition \ref{propn:wfdSless}. 
The result now follows using lemma~\ref{lem:wfdS2} and noting 
that $\bigl\lfloor\frac{\lambda}{p}\bigr\rfloor\le
\bigl\lfloor\frac{r}{p}\bigr\rfloor$ for all $\lambda \in \partn$.
\end{proof}

\begin{thm}\label{thm:glob1}
If $p>n$ then the Weyl (and the good) filtration dimension of the 
Schur algebra $S(n,r)$ is 
$$\wfd\bigl(S(n,r)\bigr) = (n-1)\left\lfloor\frac{r}{p}\right\rfloor.$$
The global dimension of $S(n,r)$ is twice this value.
\end{thm}
\begin{proof}
The previous corollary tells us that this value for $\wfd(S)$ is an upper bound for
all $p$.
So for $p>n$ we just need to give a weight in $\partn(n,r)$ whose Weyl
filtration dimension attains this bound.
We write $r=r_1p +r_0$ for $r_1$, $r_0 \in \N$ and $0\le r_0 \le p-1$.
 
Since $p> n$ we can write
$r_0=bn+a$ where $a,b\in \N$ and $0\le a\le n-1$.
Consider the weight $\mu=(r_1p +1, 1^{a-1},0^{n-a})+b(1^n)\in \partn(n,r)$.
If $a=0$ then we take $\mu=(r_1p,0^{n-1})+b(1^n)$.
The weight $\mu$ lies inside an alcove by lemma~\ref{lem:weightalc}.
Also $d(\mu)= (n-1)r_1$.
Hence $\wfd\bigl(\nabla(\mu)\bigr)=(n-1)r_1$, and so the bound is
attained.

Theorem~\ref{thm:wfdlamb} also tells us that there is a
non-zero $\Ext$ group in degree $2(n-1)r_1$. Hence the global
dimension of $S(n,r)$ is twice the Weyl filtration dimension
by corollary~\ref{cor:glob}.
\end{proof} 

\begin{thm}\label{thm:glob2}
Let $m \in \N$ then the Weyl (and the good) filtration dimension of the
Schur algebra $S(p,mp)$ is 
$$\wfd\bigl(S(p,mp)\bigr) = (p-1)m.$$
The global dimension of $S(p,mp)$ is twice this value.
\end{thm}
\begin{proof}
The weight $(mp,0,\ldots,0)\in\partn$ lies inside an alcove by
lemma \ref{lem:weightalc}.
The same argument as in the previous proof then gives us the result.
\end{proof}

The values calculated for $\wfd\bigl(L(\lambda)\bigr)$ with $\lambda$ 
inside an alcove for $n=2$ and $n=3$ agree with our previous results
for $\SL_2$ and $\SL_3$ calculated in ~\cite[sections 3 and 5]{parker1} 
This gives a new
proof for~\cite{parker1}, theorem 3.7, (for all $p$) and
\cite{parker1}, theorem 5.12, in the
cases where $p \ge 5$ and $p=3$ and $3 \mid r$.

It is still an open problem to determine what happens for weights
which are not regular. Many of the results above give upper bounds but
most of time these bounds are not sharp. Various conjectures are
presented for the value of $\wfd(S(n,r))$ in~\cite{mythesis}, section 6.5.

\section{The quantum case}\label{sect:quant}

We now show that the arguments in sections \ref{sect:wfd} and
\ref{sect:glob} generalise to the quantum case.
To do this we need the appropriate quantum versions of the results
used.
We will be using the Dipper-Donkin quantum group $q$-$\GL_n$ defined 
in~\cite{dipdonk}.  Our field $k$ remains algebraically closed but $k$
may now also have zero as well as positive characteristic.
Background information can be found in~\cite{donkbk}.
The cohomological theory of quantum groups and their $q$-Schur
algebras appears in~\cite{donkquant}.
When $q=1$ then the module category for $q$-$\GL_n$ is the same as
for $\GL_n$.
If $q$ is not a root of unity then $\Mod(q$-$\GL_n)$ is semi-simple.
We will consider the case where $q$ is a primitive $l$th root of unity with 
$l \ge 2$.

All of the structures defined in section~\ref{sect:prelim} have their
quantum analogues, which are essentially the same. The most
significant difference for us will be that $p$-alcoves and
$p$-hyperplanes will be replaced by $l$-alcoves and $l$-hyperplanes.
We need the quantum version of translation functors. These are defined
in~\cite{andpolwen}, section 8, together with the quantum version of
proposition~\ref{propn:transgood}, \cite{andpolwen}, theorem 8.3. 

All of our proofs in section \ref{sect:wfd} now carry through in the
quantum case with $p$ replaced by $l$. So the statement of
theorem~\ref{thm:wfdlamb} and \ref{thm:extiind} and their corollaries
\ref{cor:wfdnabd} and \ref{cor:wfdirr} are equally valid for the
quantum case when $l \ge h$ (even if $k$ has characteristic 0).
We also expect that the result
in the appendix carries through in the quantum case so that 
we would also have the quantum version of proposition \ref{propn:nabnab} and
theorem \ref{thm:injproj}.

We now consider the quantised Schur algebra, $S_q(n,r)$.
This can be constructed in the same way as in the classical case.
Take the saturated subset of dominant weights
$\Pi=\partn(n,r)$, then the quasi-hereditary algebra
$S(\Pi)$ is isomorphic to $S_q(n,r)$, the quantised Schur algebra.
Moreover we have the same ordering -- namely the $\up$-ordering
defined using the action of the affine Weyl group.
The proofs in section \ref{sect:glob} work equally well in the 
quantum case with $p$ replaced by
$l$ where $q$ is an $l$th root of unity with $l \ge 2$. So we get an upper bound for
$\wfd\bigl(S_q(n,r)\bigr)$ of $(n-1)\bigl\lfloor \frac{r}{l} \bigr\rfloor$.
Together with the quantum version of the results of  section
\ref{sect:wfd} we may now deduce the following theorem.
\begin{thm}
If $q$ is a primitive $l$th root of unity with $l>n$ then the Weyl 
(and the good) filtration dimension of the quantised 
Schur algebra $S_q(n,r)$ is 
$$\wfd\bigl(S_q(n,r)\bigr) =
(n-1)\left\lfloor\frac{r}{l}\right\rfloor.$$
Suppose $l \ge 2$ and let $m \in \N$. Then we have
$$\wfd\bigl(S_q(l,ml)\bigr) = (l-1)m.$$
In both these case the global dimension of $S_q(n,r)$ and $S_q(l,ml)$
is twice its Weyl filtration dimension.
\end{thm}
Again, this result is dependent only on $l$ and not on the
characteristic of the field $k$.

\section{Category $\OO$}
There are analogous situations in Category $\OO$ defined by
Bern\v ste\u\i n, Gel'fand and Gel'fand,~\cite{BGG2}.
Category $\OO$ is known to be a highest weight category
(see \cite{klukoe}, section 4.1 for a basic introduction) 
so we can apply the general theory of section~\ref{sect:qha}. 
We use the setup of~\cite{carlin}, although  
note that~\cite{carlin} uses the terminology `$p$-filtration' for what
we have defined to be a Weyl filtration.
There $\g$ is a complex, semi-simple Lie algebra with Cartan
subalgebra $\h$ and Weyl group $W$. We denote the longest element of
$W$ by $w_0$.
The standard modules for $\OO$ are the well-known Verma modules,
denoted $M(\lambda)$ for $\lambda \in \h^*$. 
We also have that $[M(\mu):L(\lambda)]\ne 0$
if and only if there are positive roots $\gamma_1,\ldots, \gamma_m$
such that there is a chain of inequalities $\mu \ge s_{\gamma_1}(\mu)
\ge \cdots \ge s_{\gamma_m}\cdots s_{\gamma_1}(\mu) = \lambda$ (\cite{BGG1}).
We may use \cite{carlin}, proposition 3.7, theorem
3.8 and theorem 4.6 to deduce that $\gfd\bigl(M(w \cdot
\lambda)\bigr)=\gfd\bigl(L(w\cdot \lambda)\bigr)=l(w_0)-l(w)$  and
$\proj\bigl(M(w\cdot \lambda)\bigr) = l(w)$ for
$\lambda$ an integral weight inside the dominant Weyl chamber. We may deduce 
that $\proj\bigl(L(w \cdot \lambda)\bigr) \le
2l(w_0)- l(w)$. These last two statements are consistent with \cite{BGG2}, statements
1 and 2. We also have translation functors and the
analogue of proposition~\ref{propn:transgood} and hence the
corollary~\ref{cor:transses}.
Unfortunately the analogue of lemma~\ref{lem:wfdmu} may no longer be true. 
Our argument does not work in this situation and indeed already fails
for type $A_2$. However, in~\cite{BGG2}, remark in \S7, it is stated
that $\Ext^{2l(w_0)}\bigl(L(\lambda), L(\lambda)\bigr) \cong \C$.
So there is strong evidence to suggest that
$\Ext^{2l(w_0)-l(w)-l(v)}\bigl(L(\wlamb),L(\vlamb)\bigr) \cong \C$ for $v$, $w \in
W$.
The results of~\cite{BGG2}, \S7, are already enough to deduce that
the global dimension of $\OO$ is $2l(w_0)$.


\begin{thebibliography}{10}

\bibitem{ander1}
H.~H. Andersen, \emph{The strong linkage principle}, J. reine angew. Math.
  \textbf{315} (1980), 53--59.

\bibitem{andpolwen}
H.~H. Andersen, P.~Polo, and K.~X. Wen, \emph{Representations of quantum
  algebras}, Invent. Math. \textbf{104} (1991), no.~1, 1--59.

\bibitem{benson}
D.~J. Benson, \emph{{R}epresentations and {C}ohomology {I}}, Cambridge Studies
  in Advanced Mathematics, no.~30, Cambridge University Press, 1995.

\bibitem{BGG1}
I.~N. {Bern\v ste\u\i n}, I.~M. Gel'fand, and S.~I. Gel'fand, \emph{Structure
  of representations generated by highest weight}, Funct. Anal. and Appl.
  \textbf{5} (1971), 1--8.

\bibitem{BGG2}
\bysame, \emph{A category of {${\mathfrak g}$}--modules}, Funct. Anal. and
  Appl. \textbf{10} (1976), 87--92.

\bibitem{carlin}
K.~J. Carlin, \emph{Extensions of {V}erma modules}, Trans. Amer. Math. Soc.
  \textbf{294} (1986), no.~1, 29--43.

\bibitem{dipdonk}
R.~Dipper and S.~Donkin, \emph{Quantum {$\mathrm{GL}_n$}}, Proc. London Math.
  Soc. (3) \textbf{63} (1991), 165--211.

\bibitem{donkrat}
S.~Donkin, \emph{{Rational Representations of Algebraic Groups: Tensor Products
  and Filtrations}}, Lecture Notes in Mathematics, vol. 1140, Springer--Verlag,
  Berlin/Heidelberg/New York, 1985.

\bibitem{donk1}
\bysame, \emph{On {S}chur algebras and related algebras {I}}, J. Algebra
  \textbf{104} (1986), 310--328.

\bibitem{donktilt}
\bysame, \emph{On tilting modules for algebraic groups}, Math. Z. \textbf{212}
  (1993), 39--60.

\bibitem{donkquant}
\bysame, \emph{Standard homological properties for quantum {$\mathrm{GL}_n$}},
  J. Algebra \textbf{181} (1996), 235--266.

\bibitem{donkbk}
\bysame, \emph{The {$q$}--{S}chur {A}lgebra}, London Math. Soc. Lecture Note
  Ser., vol. 253, Cambridge University Press, Cambridge, 1998.

\bibitem{dotynak}
S.~R. Doty and D.~K. Nakano, \emph{Semi--simplicity of {S}chur algebras}, Math.
  Proc. Cambridge Philos. Soc. \textbf{124} (1998), 15--20.

\bibitem{erdnak}
K.~Erdmann and D.~K. Nakano, \emph{Representation type of $q$--{S}chur
  algebras}, Trans. Amer. Math. Soc. \textbf{353} (2001), no.~12, 4729--4756.

\bibitem{fripar}
E.~M. Friedlander and B.~J. Parshall, \emph{Cohomology of {L}ie algebras and
  algebraic groups}, {A}mer. {J}. {M}ath. \textbf{108} (1986), 235--253.

\bibitem{humph}
J.~E. Humphreys, \emph{{L}inear {A}lgebraic {G}roups}, Graduate Texts in
  Mathematics, vol.~21, Springer--Verlag, Berlin/Heidelberg/New York, 1975.

\bibitem{humph2}
\bysame, \emph{{R}eflection {G}roups and {C}oxeter {G}roups}, Cambridge Studies
  in Advanced Mathematics, no.~30, Cambridge University Press, 1990.

\bibitem{jantz}
J.~C. Jantzen, \emph{{R}epresentations of {A}lgebraic {G}roups}, Pure Appl.
  Math., vol. 131, Academic Press, San Diego, 1987.

\bibitem{klukoe}
M.~Klucznik and S.~K\"onig, \emph{{Characteristic Tilting Modules over
  Quasi--hereditary Algebras}}, unpublished notes, 1999.

\bibitem{mathieu}
O.~Mathieu, \emph{Filtrations of {$G$}--modules}, Ann. Sci. \'Ecole Norm. Sup.
  (4) \textbf{23} (1990), no.~4, 625--644.

\bibitem{parker1}
A.~E. Parker, \emph{The global dimension of {S}chur algebras for
  {$\mathrm{GL}_2$} and {$\mathrm{GL}_3$}}, J. Algebra \textbf{241} (2001),
  340--378.

\bibitem{mythesis}
\bysame, \emph{On the global dimension of {S}chur algebras and related
  algebras}, Ph.D. thesis, University of London, 2001.

\bibitem{springer}
T.~A. Springer, \emph{{L}inear {A}lgebraic {G}roups}, Progress in Mathematics,
  vol.~9, Birkh\"auser, Boston/Basel/Stuttgart, 1981.

\bibitem{totaro}
B.~Totaro, \emph{Projective resolutions of representations of
  {${\mathrm{GL}}(n)$}}, J. reine angew. Math. \textbf{482} (1997), 1--13.

\bibitem{verma}
D.~N. Verma, \emph{The {r\^ole} of affine {W}eyl groups}, {L}ie Groups and
  their representations (I.M. Gel'fand, ed.), 1975, pp.~653--705.

\bibitem{wang1}
Jian-pan Wang, \emph{Sheaf cohomology of {$G/B$} and tensor products of {W}eyl
  modules}, J. Algebra \textbf{77} (1982), 162--185.

\end{thebibliography}

\providecommand{\bysame}{\leavevmode\hbox to3em{\hrulefill}\thinspace}
\providecommand{\MR}{\relax\ifhmode\unskip\space\fi MR }
\providecommand{\MRhref}[2]{%
  \href{http://www.ams.org/mathscinet-getitem?mr=#1}{#2}
}
\providecommand{\href}[2]{#2}

\end{document}